# A Multi-stage Stochastic Programming Approach for Pre-positioning of Relief Supplies.

## Abstract ID: 581648

Oluwasegun Olanrewaju, Shaolong Hu, Zhijie Sasha Dong
Ingram School of Engineering, Texas State University
San Macros, TX 78666

## Abstract
Pre-positioning of relief supplies is an important aspect of disaster operations management that aims at decreasing the response time by advancing procurement and storage of needed supplies. In this paper we consider commodity life-time period with the related costs in keeping the commodity and removing it from the storage when it is close to expiration (i.e. holding and removal cost). We develop a multi-stage stochastic programming model for pre-positioning of relief supplies, which provides relief agencies insight on how to have dynamically control inventories due to uncertain demands when disasters (e.g., hurricane or earthquake) occur. We also present a case study based on mainland United States to illustrate the proposed model as well as provide managerial insights to relief agencies.

## Keywords
Emergency logistics, pre-positioning, inventory control, perishable goods, multi-stage stochastic programming

## 1. Introduction
Disaster operations management is an important aspect which contributes to the improvement in readiness to disaster, reducing injuries, fatalities, and damages caused by disasters, and to ease recovery. The focus of this paper is on the preparedness stage, precisely on pre-positioning of relief supplies before the occurrence of a disaster in the United States. In the wake of natural disaster there are always many people affected. To reduce fatalities, the needed assistance must be rendered timely after the occurrence of the disaster. Therefore, there is always a tremendous amount of demand for various basic relief supplies such as water, food and medical kits. This brings about the need to have a good plan in place before disasters occur. Pre-positioning of relief supplies is a key activity in disaster operations management, which helps preparedness for natural disasters by advancing procurement of needed supplies, thereby decreasing the response time. However, this is challenging, because pre-positioning requires high investment (e.g., procurement and holding costs) at various locations, due to a high level of uncertainty in the timing and location of the next disaster. In addition, limited life time of relief items is another major problem [1]. For example, when the items in the warehouse have not been used due to nonoccurrence of disaster, and it is close to its expiration date, the relief agency must donate or dispose the items. This study focusses on pre-positioning of relief supplies by considering uncertain and dynamic demand for products which have expiration dates. It also provides relief agencies (e.g. FEMA) managerial insights about dynamic control of inventory over each scenario and dealing with relief supplies which will expire. The paper is organized as follows. Section 2 entails review of relevant literature. A multi-stage stochastic programming model is formulated in Section 3. Section 4 illustrates the numerical analysis with data setting and numerical results. Conclusion and future research directions are summarized in Section 5.

## 2. Literature review
Due to the frequent occurrence of natural disasters in the world, pre-positioning of relief supplies has been a major area that researchers have dived into. [2] developed a new integrated model that determines the optimum location-allocation and distribution plan, coupled with the best ordering policy of renewing the stocked perishable commodities at the disaster phase. They were concerned with the periodic ordering policy of commodity while considering a fixed lifetime for the commodity. [3] proposed a multi-stage stochastic programming model for relief distribution by considering the state of road network and multiple types of vehicles are considered. [4] integrated supplier selection into the decision making of pre-positioning relief supplies, and formulated the problem as a two-stage stochastic programming model with a consideration of failure risks under disasters. [1] discussed how pre-positioning of relief inventory requires a high investment and proposed a dynamic model that entails the delivery process of ready-to-use therapeutic food items during the immediate response phase of a disaster. They also analyzed the performance of different preparedness scenarios. [5] proposed a stochastic linear mixed-integer programming model for integrated decisions in the preparedness and response stages in pre- and post-disaster operation respectively, while considering key areas such as facility and stock pre-positioning, evacuation planning and relief vehicle planning. [6] presented a forecast driven dynamic model for pre-positioning of relief items in preparation for a foreseen hurricane. The model uses forecast advisories which helps update the information to determine the amount and location of units to be pre-positioned and re-prepositioned. [7] developed an emergency response planning tool that determines the location and quantities of various types of emergency supplies to be prepositioned. They presented a two-stage stochastic mixed integer program (SMIP) that helps in the pre-positioning strategy for hurricanes. [8] discussed about how FEMA operates as well as the importance of distributions recovery centers and their roles during disasters. All these studies are important in the literature on pre-positioning of relief supplies. However, the major limitation of their work is that they assumed that inventory is fixed, and they do not consider that relief items will expire. Compared to literature in pre-positioning of relief supplies, this paper has two distinctive contributions. Firstly, we consider the life-time of a commodity type and the cost associated with removing the item when it is close to expiration. This makes the model closer to real situations. Secondly, we formulate a multi-stage stochastic programming model which enables dynamic and stochastic control over the inventory.

## 3. Modeling

A multi-stage stochastic programming model is created such that the uncertainty of demand for relief supplies is presented through a scenario tree. The objective of this model is to determine sequencial decisions on the amount of relief items to be procured and the balancing of the inventories. Fig. 1 shows the scenario tree which is planned over a time horizon, representing a discrete set of scenarios. In each statge, there could be several different sceanrios, and only one scenario would be realized. The lifetime period corresponds to the time between each stage, which can be days, months, and years. The first scenario, which is the first stage, is denoted by a root scenario. Any scernario occurs after a specific scenario is called the child scenrio of this scenario, and this specific sceanrio is the parent scenario of the subsequent sceairos. For exampel, in Fig. 1, scenario 1 is the parent scenario of scenario 2, scenario 2 is the parent scenario of scenario 5, and scenario 5 is the parent scenario of scenario 9. The sum of probabilities of each scenario at a given stage is equal to one because the probability of occurrence of each child state of a given parent scenario has already happened in the previous stage.

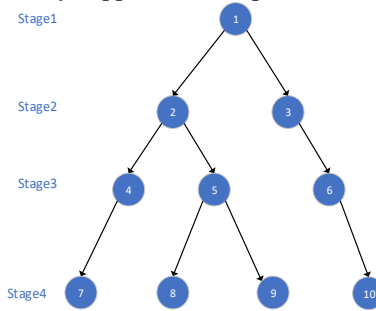

Fig. 1: Example of a scenario tree

Table 1. Variation of inventory for relief supplies with lifetime and scenario

| Stage | Scenario | $t=1$ (12 months) | $t=2$ (8 months) | $t=3$ (4 months) | $t=4$ (0 month) |
|---|---|---|---|---|---|
| 1 | 1 | $x(1) \rightarrow h(1,1)$ | 0 | 0 | 0 |
| 2 | 2 | $x(2) \rightarrow h(2,1)$ | $h(2,2)$ | 0 | 0 |
| 3 | 5 | $x(5) \rightarrow h(5,1)$ | $h(5,2)$ | $h(5,3)$ | 0 |
| 4 | 9 | 0 | $h(9,2)$ | $h(9,3)$ | $h(9,4)$ |

The lifetime of the commodity is discretized, and $T^c$ is introduced to represent the set of remainig lifetime period for commodity type $c$, where $t \in T^c$. $T^c$ represents a general setting because relief items usually have different lifetimes. For example, if we assume four remaining lifetime periods for commidy type 1, $T^1 = \{1,2,3,4\}$, it means that these are items with 12-months total lifetime and the lapse of time between remaining lifetime periods represents a 4-months life interval. Then $t=1$ represents 12 remaining months, $t=2$ represents 8 remaining months, $t=3$ represents remaining 4 months, and $t=4$ represents 0 months remaining (i.e., expired). Table1 shows details about the variation of inventory over a lifetime for scenario tree 1-2-5-9, where $x(s)$ represents the procurement quantity at scenario $s$, $x(s,t)$ represents the procurement quantity at scenario $s$ at lifetime period $t$, $d(s,t)$ represents the commodity used or demanded at lifetime period $t$ at scenario $s$, and $h(s,t)$ represents the inventory at lifetime period $t$ at scenario $s$, respectively. We assume that we only procure the relief supplies with a maximum lifetime period (i.e., 12 months in the previous example). Therefore, inventory at remaining lifetime period 1 will be equal to procurement quantity. Inventories at other remaining lifetime periods depend on $x(s,t)$ and $d(s,t)$. The decisions made on $x(s)$ and $x(s,t)$ depend on the probabilities and demands of subsequent scenarios. For example, the inventory at scenario 2 when remaining lifetime period is 2 ($h(2,2)$) depends on the procurement quantity and the commodity used at the same scenario and lifetime (i.e., $h(2,2) = h(1,1) - x(2,2) - d(2,2)$, $q(1) > 0$). It is worth mentioning that for those relief supplies which have been procured in scenario 1 will be disposed by the end of stage 4 and scenario 9 (i.e., when $h(9,4)>0$) because they are expired. It is the same with scenarios 7, 8, and 10.

The following notations are used for the multi-stage stochastic programming model:

The main set are:

| | | | |
|---|---|---|---|
| $C$ | set of commodities | $S$ | set of discrete scenarios |
| $T^c$ | set of remaining lifetime periods for type $c$ items | | |

The Indices of sets are:

| | | | |
|---|---|---|---|
| $c$ | type of commodity $c \in C$ | $s$ | scenario type $s \in S$ |
| $k$ | parent scenario $k \in K^s$ | $i,j$ | specific location $i, j \in N$ |
| $t$ | remaining lifetime period for type $c$ commodity $t \in T^c$ | | |

The deterministic parameters are:

| | | | |
|---|---|---|---|
| $M_i$ | overall capacity of facility at location $i$ | $b^c$ | unit space required for commodity $c$ |
| $q^c$ | unit acquisition cost of commodity $c$ | $r^c$ | unit removal cost of commodity $c$ |
| $o^c_{i,j}$ | unit cost of transporting commodity $c$ across link $(i,j)$ | $u^c$ | holding cost of unused commodity $c$ |
| $v^c$ | penalty cost associated with shortage of commodity $c$ | | |

The decision variables are:

| | | | |
|---|---|---|---|
| $x^{cs}_i$ | procurement quantity of commodity $c$ at location $i$ in scenario $s$ | $y^{cs}_{ijt}$ | amount of commodity $c$ shipped across the link $(i,j)$ within remaining lifetime period $t$ in scenario $s$ |
| $g^{cs}_j$ | shortage of commodity $c$ at location $j$ in scenario $s$ | $h^{cs}_{it}$ | current inventory of commodity $c$ at location $i$ within remaining lifetime period $t$ in scenario $s$ |

The extensive formulation of the model is presented below:

$$\min \sum_s p_s \left[ \sum_{c,i} q^c x^{cs}_i + \sum_{c,i,j,t} o^c_{ij} y^{cs}_{ijt} + \sum_{t<|T^c|,c,i} u^c h^{cs}_{it} + \sum_{c,j} v^c g^{cs}_j + \sum_{c,i} r^c h^{cs}_{i|T^c|} \right] \quad (1)$$

$$h^{c1}_{it} = 0, \quad t \geq 2, \forall i,c \quad (2)$$

$$h^{cs}_{i1} = x^{cs}_i, \quad \forall i,c,s \quad (3)$$

$$h^{cs}_{it} = h^{ck}_{it-1} - \sum_j y^{cs}_{ijt}, \quad s,t \geq 2, \forall i,k,c \quad (4)$$

$$y^{cs}_{ij1} = 0, \quad \forall c,i,j,s \quad (5)$$

$$\sum_{ct} b^c h^{cs}_{it} \leq M_i, \quad \forall s,i \quad (6)$$

$$g^{cs}_j = d^{cs}_j - \sum_{i,t} y^{cs}_{ijt}, \quad \forall j,c,s \quad (7)$$

$$x^{cs}_i, y^{cs}_{ijt}, g^{cs}_j, h^{cs}_{it} \geq 0, \quad \forall i,j,c,t,s \quad (8)$$

The objective function (1) minimizes the total expected cost over all scenarios resulting from the selection of the pre-positioning locations and facility sizes, the commodity acquisition and stocking decisions, the shipments of the supplies to the demand points, unmet demand penalties and holding costs for unused material. Constraint (2) limits the the current inventory of the first scenario for all location to be zero. Constraint (3) restricts the current inventory of scenario $s$ at the first stage to be the same as the quantity of the procured commodity at location $i$. Constraint (4) states that the current inventory at a given stage is the difference between the inventory of the previous stage and amount of commodity shipped out. Constraint (5) limits the amount of commodities shipped out at the first stage to zero. Constraint (6) restricts the space occupied by the stocked commodities not exceed the facility capacity at location $i$. Constraint (7) shows how shortage of commodities is calculated. Constraint (8) shows sign of decision variables.

## 4. Numerical Analysis
### 4.1 Data Setting
The focus of this research is on pre-positioning of relief supplies in preparedness for hurricane, flood and earthquake threat in the United State. We study 49 states in mainland United States as demand locations. The distance between states are calculated using the most populous city in each state. Two commodities are considered: food and water. The unit of water and food is assumed to be 1,000 gallons and 1,000 meals ready to eat (MREs), respectively. For each commodity type a remaining life-time period of four is considered. The facility locations considered are Texas (1.6 million sq. ft), California (110,000 sq. ft), Georgia (407,000sq. ft) and Maryland (68,023sq. ft). The facilities are already provided by Federal Emergency Management Agency (FEMA) [12]. The values of acquisition cost $q^c$, space estimate $b^c$, transport cost $o_{i,j}^c$, penalty cost $v^c$, removal cost $r^c$ and holding cost $u^c$, are provided in Table 2. The cost of transporting a commodity between two locations is derived by multiplying the unit transport cost assigned to the commodity type and the distances between the location. The following parameter follows [7]. The holding cost is assumed to be 25% of the acquisition cost of commodities. The penalty cost for unmet demand is assumed to be 10 times the acquisition cost of commodities. The removal cost is assumed to be 40% of the acquisition cost of commodities.

**Table 2. Unit procurement price, transport, holding, removal & penalty cost and storage volume occupied**

|  | $b^c$(ft$^3$/unit) | $q^c$($/unit) | $u^c$($/unit) | $v^c$($/unit) | $r^c$($/unit) | $o_{i,j}^c$($/unit-mile) |
|---|---|---|---|---|---|---|
| Water (1000 gals) | 144.6 | 647.7 | 161.925 | 6477 | 259.08 | 0.3 |
| Food(1000MREs) | 83.33 | 5420 | 1355 | 54200 | 2168 | 0.04 |

For the emergency demand data of each location, we used the top 10 ranked location [11] with the highest occurrence of hurricane, earthquake and flood. Note that we do not consider the specific intensity for hurricane, earthquake and flood. To simplify data estimation, we assume three impact levels for each disaster, such as low, medium and high. The demand for hurricane and flood is generated from uniform distribution(U) by assigning U [100,200], U [200,400] and U [900,1000] to low, medium and high impact hurricanes and flood, respectively. The demand for earthquake is generated by assigning U [1000,2000], U [2000,4000] and U [9000,10000] to low, medium and high impact respectively [11]. The demands are computed based on the information on hurricane, flood *and* earthquake together. We assume the same demand for both relief items. The model is solved using CPLEX Concert Technology (IBM ILOG CPLEX) in Microsoft Visual Studio as an integrated development environment (IDE) on a laptop with Intel Core i7-700 @3.60GHz and 8GB RAM. The optimal solution was obtained with 139.53 seconds.

### 4.2 Numerical Results
In this section, sensitivity analysis is carried out in order to examine how different parameters in the proposed model will cause changes and thus affect decision making about pre-positioning relief supplies as well as all the costs. The unit penalty cost $v^c$, unit removal cost $r^c$ and unit holding cost $u^c$ are the parameters that will be studied. We are particularly interested in the effects of varying parameters on the economic cost (which is the total acquisition cost $Q$ total transport cost $O$, total holding cost $U$ and total removal cost $R$) and the total penalty cost $V$.

**(1) Effects of holding cost**

**Table 3. Sensitivity of costs to unit holding cost and unit penalty cost**

|  | Unit holding cost | | | | Unit penalty cost | | | |
|---|---|---|---|---|---|---|---|---|
|  | -50% | -25% | +25% | +50% | -50% | -25% | +25% | +50% |
| $Q$ | +3.08 | +2.17 | -0.31 | -2.48 | -64.60 | -6.20 | +3.19 | +3.67 |
| $O$ | +2.36 | +1.70 | -0.30 | -2.47 | -53.77 | -5.68 | +2.54 | +3.12 |
| $U$ | -48.12 | -23.20 | +24.57 | +46.41 | -66.32 | -5.55 | +3.85 | +5.51 |
| $V$ | -3.86 | -2.62 | +0.42 | +3.16 | +7.89 | -18.99 | +20.06 | +42.35 |
| $R$ | +5.12 | +2.64 | -0.51 | -2.12 | -77.85 | -4.83 | +5.17 | +9.36 |
| *Total* | -9.34 | -4.64 | +4.56 | +9.07 | -32.36 | -11.88 | +11.19 | +22.15 |

Because of uncertainty of disaster occurrence, the longer time the prepositioned relief items spend in the warehouse, the larger the holding cost will be. When disasters occur earlier than expected or much later, the cost associated with keeping it safe in the store/warehouse, which is known as holding cost, will be affected. Thus, we modify the unit holding cost from -50% of the original inputs to +50% and the results obtained is summarized in Table 3. Holding cost is directly proportional to total penalty cost. But total acquisition and removal cost decreases with increase in unit holding cost, because the procurement quantity decreases. Therefore, the relief agencies won't purchase new commodities due to increment in the holding cost. This in return leads to less relief item that will be dispose due to closeness to expiration. Figure 2 shows the behavior of total penalty cost and economic cost when the unit holding cost is being modified. The increase in unit holding cost directly affect the economic costs. The total penalty cost tends to increase when the unit holding cost increases because an increment in the holding cost indicates the relief items are still in the warehouse and not been shipped to a place where it is needed due to nonoccurrence of disaster.

**(2) Effects of penalty cost**

We modify the unit penalty cost from -50% of the original Input to +50% and the effects of this change on different costs terms are illustrated in Table 3. Increase in unit penalty cost leads to the increase of the acquisition, removal and holding cost. The relief agencies tend to purchase more relief items when penalty associated with shortage of relief supplies is high. This then leads to the increase of relief items in store and if not use leads to increase of relief items to dispose when it is almost expired. Figure 3 shows how economic cost and total penalty cost changes when penalty cost changes. When the unit penalty cost is at its lowest, the shortage of commodity at different locations increases. Relief agencies tend not to feel the urgency to ship commodity to a location if the consequence is negligible compared to when the penalty associated with not shipping the commodity is high.

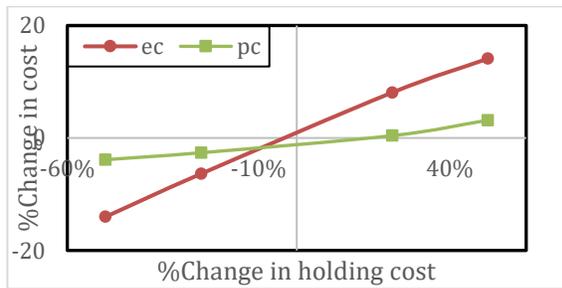

**Figure 2. Observation on economic cost and total penalty cost when varying the holding cost.**

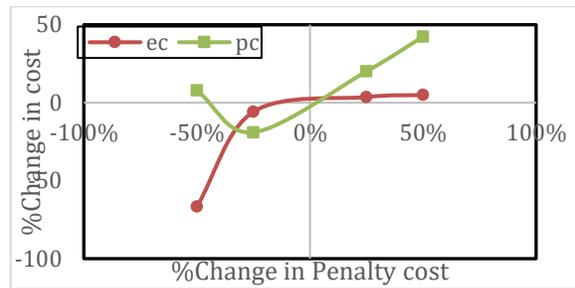

**Figure 3. Observation on economic cost and total penalty cost when varying the penalty cost.**

**(3) Effects of removal cost**

**Table 4. Sensitivity of costs to removal cost**

|       | -50%   | -25%   | +25%   | +50%   | +75%   | +100%  |
|-------|--------|--------|--------|--------|--------|--------|
| Q     | +1.73  | +1.51  | -0.27  | -0.31  | -0.62  | -1.11  |
| O     | +1.55  | +0.81  | -0.03  | -0.29  | -0.82  | -1.12  |
| U     | +1.97  | +1.66  | -0.31  | -0.34  | -0.68  | -1.23  |
| V     | -2.41  | -1.91  | +0.37  | +0.42  | +0.86  | +1.60  |
| R     | -47.38 | -23.23 | +24.42 | +49.23 | +73.27 | +96.42 |
| Total | -3.19  | -1.58  | +1.55  | +3.10  | +4.64  | +6.18  |

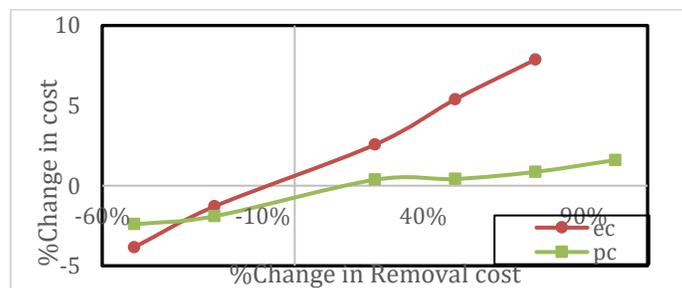

**Figure 4. Observation on economic cost and total penalty cost when varying the removal cost.**

Removal cost is the cost associated with removing the relief item from the store or warehouse when its expiration date gets close. There is tendency for it to vary since it depends on whether disasters occur or not, and when it does, whether the items are put into use within their life-time period. To verify the effect, we modify the unit removal cost from -50% of the original inputs to +100% and the results obtained are summarized in Table 4. The total acquisition cost and holding cost also decreases as the unit removal cost increases. This is due to the decrease in procurement quantity with the aim of reducing the amount of commodity close to expiration. Figure 4 shows the behavior of total penalty cost, economic cost when the unit removal cost is being modified. With the decrease in removal cost, the economic cost and the total penalty cost decreases. When disasters occur, the relief agent tends to use up commodities that have a shorter life-time period, which helps in prolonging the time the remaining commodities would be removed while waiting on disaster occurrence. When unit removal cost increases the total penalty cost tends to increase, because the procurement quantity reduces.

## 5. Conclusions

This paper addresses the inventory management of two basic relief supplies during the preparedness stage of disaster management by considering the limited life period of supplies. A multi-stage stochastic programing model is formulated with the objective of minimizing the total cost, which consists of the acquisition cost, space estimate, transport cost, penalty cost, removal cost, and holding cost. The major contribution of this paper is that it takes into consideration, uncertain demand and life time of products, thereby giving relief agency insights in dynamic control of inventory over each scenario and dealing with relief items which will expire.

Sensitivity analysis of experiment comparison has provided the relief agencies insights on how they can manage the cost associated with the preposition of relief supplies while having dynamic control over the inventory. Firstly, avoiding shortage of a commodity is a good cost-effective process because of the penalty associated with not having the commodity. Secondly, having a better and more efficient way of disposing commodity when they have low lifetime period is also an aspect the relief agencies should consider. Thirdly, relief agencies should have a better plan in place for handling relief items, because the more time it spends in the warehouse the larger the cost of holding it.

Future research areas include bringing in suppliers as an alternative of removing relief items when its expiration is close instead of disposing items. To be more specific, commodities close to expiration are returned to suppliers to sell, and in the event that disasters occur, the commodities will be replaced by suppliers.